\documentclass[11pt]{amsart}
\usepackage{amsmath,amssymb,stmaryrd,mathrsfs}
\newtheorem{theorem}{Theorem}
\newtheorem{corollary}[theorem]{Corollary}
\newtheorem{lemma}[theorem]{Lemma}
\newtheorem{proposition}[theorem]{Proposition}
\newcommand{\scal}{\text{\rm scal}}
\newcommand{\Ric}{\text{\rm Ric}{}}
\newcommand{\tracefreeRic}{\overset{\text{\rm o}}{\text{\rm Ric}}{}}

\begin{document}

\title{Ancient solutions to the Ricci flow with pinched curvature}
\author{Simon Brendle, Gerhard Huisken, and Carlo Sinestrari}
\address{Department of Mathematics \\ Stanford University \\ Stanford, CA 94305}
\address{Max-Planck Institut f\"ur Gravitationsphysik, Am M\"uhlenberg 1, 14476 Golm, Germany}
\address{Dipartimento di Matematica, Universit\`a di Roma "Tor Vergata", Via della Ricerca Scientifica, 00133 Roma, Italy}
\maketitle 

\section{Introduction}

In this note, we study ancient solutions to the Ricci flow on compact manifolds. Recall that a one-parameter family of metrics $g(t)$ on a compact manifold $M$ evolves by the Ricci flow if 
\[\frac{\partial}{\partial t} g = -2 \, \Ric_{g(t)}.\] 
A solution to the Ricci flow is called ancient if it is defined on a time interval $(-\infty,T)$. Ancient solutions typically arise in the study of singularities to the Ricci flow (see e.g. \cite{Hamilton3}, \cite{Hamilton-survey}, \cite{Perelman1}, \cite{Perelman2}).

P.~Daskalopoulos, R.~Hamilton, and N.~\v Se\v sum \cite{Daskalopoulos-Hamilton-Sesum2} have recently obtained a complete classification of all ancient solutions to the Ricci flow in dimension $2$. (See also \cite{Daskalopoulos-Hamilton-Sesum1}, where the analogous question for the curve shortening flow is studied.) V.~Fateev \cite{Fateev} has constructed an interesting example of an ancient solution in dimension $3$. L.~Ni \cite{Ni} showed that any ancient solution to the Ricci flow which is of Type I, $\kappa$-noncollapsed, and has positive curvature operator has constant sectional curvature.

In this note, we show that any ancient solution to the Ricci flow in dimension $n \geq 3$ which satisfies a suitable curvature pinching condition must have constant sectional curvature. In dimension $3$, we require a uniform lower bound for the Ricci tensor:

\begin{theorem}
\label{ancient.solutions.dim.3}
Let $M$ be a compact three-manifold, and let $g(t)$, $t \in (-\infty,0)$, be an ancient solution to the Ricci flow on $M$. Moreover, suppose that there exists a uniform constant $\rho > 0$ such that 
\[\Ric_{g(t)} \geq \rho \, \scal_{g(t)} \, g(t) \geq 0\] 
for all $t \in (-\infty,0)$. Then the manifold $(M,g(t))$ has constant sectional curvature for each $t \in (-\infty,0)$.
\end{theorem}

Fateev's example shows that the pinching condition for the Ricci tensor cannot be removed. The proof of Theorem \ref{ancient.solutions.dim.3} relies on a new interior estimate for the Ricci flow in dimension $3$ (cf. Proposition \ref{interior.estimate.dim.3} below). The proof of this estimate relies on the maximum principle, and will be presented in Section \ref{dim.3}.

In dimension $n \geq 4$, we prove the following result:

\begin{theorem}
\label{ancient.solutions.higher.dim}
Let $M$ be a compact manifold of dimension $n \geq 4$, and let $g(t)$, $t \in (-\infty,0)$, be an ancient solution to the Ricci flow on $M$. Moreover, suppose that there exists a uniform constant $\rho > 0$ with the following property: for each $t \in (-\infty,0)$, the curvature tensor of $(M,g(t))$ satisfies 
\begin{align*} 
&R_{g(t)}(e_1,e_3,e_1,e_3) + \lambda^2 \, R_{g(t)}(e_1,e_4,e_1,e_4) \\ 
&+ R_{g(t)}(e_2,e_3,e_2,e_3) + \lambda^2 \, R_{g(t)}(e_2,e_4,e_2,e_4) \\ 
&- 2\lambda \, R_{g(t)}(e_1,e_2,e_3,e_4) \geq \rho \, \scal_{g(t)} \geq 0 
\end{align*}
for all orthonormal four-frames $\{e_1,e_2,e_3,e_4\}$ and all $\lambda \in [0,1]$. Then the manifold $(M,g(t))$ has constant sectional curvature for each $t \in (-\infty,0)$.
\end{theorem}

Theorem \ref{ancient.solutions.higher.dim} again follows from pointwise curvature estimates which are established using the maximum principle (see Corollary \ref{interior.estimate.higher.dim} below). In dimension $n \geq 4$, the evolution equation for the curvature tensor is much more complicated, and our estimates are not as explicit as in the three-dimensional case. In order to handle the higher dimensional case, we use the invariant curvature conditions introduced in \cite{Brendle1} and \cite{Brendle-Schoen1}. These ideas also play a key role in the proof of the Differentiable Sphere Theorem (cf. \cite{Brendle-Schoen1}, \cite{Brendle-Schoen-survey}).

\section{Proof of Theorem \ref{ancient.solutions.dim.3}}

\label{dim.3}

\begin{proposition}
\label{interior.estimate.dim.3}
Let $M$ be a compact three-manifold, and let $g(t)$, $t \in [0,T)$, be a solution to the Ricci flow on $M$. Moreover, suppose that there exists a uniform constant $\rho \in (0,1)$ such that 
\[\Ric_{g(t)} \geq \rho \, \scal_{g(t)} \, g(t) \geq 0\] 
for each $t \in [0,T)$. Then, for each $t \in (0,T)$, the curvature tensor of $(M,g(t))$ satisfies the pointwise estimate 
\[|\tracefreeRic_{g(t)}|^2 \leq \Big ( \frac{3}{2t} \Big )^\sigma \, \scal_{g(t)}^{2-\sigma},\] 
where $\sigma = \rho^2$.
\end{proposition}

\textbf{Proof.} 
The assertion is trivial if $(M,g(0))$ is Ricci flat. Hence, it suffices to consider the case that $(M,g(0))$ is not Ricci flat. By the maximum principle, the manifold $(M,g(t))$ has strictly positive scalar curvature for all $t \in (0,T)$. 

We next define a function $f: M \times (0,T) \to \mathbb{R}$ by 
\[f = \scal^{\sigma-2} \, |\tracefreeRic|^2,\] 
where $\sigma = \rho^2$. It is easy to see that $f \leq \scal^\sigma$. Moreover, it follows from Lemma 10.5 in \cite{Hamilton1} that 
\[\frac{\partial}{\partial t} f \leq \Delta f + \frac{2(1-\sigma)}{\scal} \, \partial_k \scal \, \partial^k f + 2 \, \scal^{\sigma-3} \, \Big [ \sigma \, |\Ric|^2 \, |\tracefreeRic|^2 - 2P \Big ],\] 
where $P$ is a polynomial expression in the eigenvalues of the Ricci tensor. By assumption, we have $\Ric \geq \rho \, \scal \, g$. Hence, it follows from Lemma 10.7 in \cite{Hamilton1} that 
\[P \geq \sigma \, |\Ric|^2 \, |\tracefreeRic|^2.\] 
This implies 
\begin{align*} 
2P - \sigma \, |\Ric|^2 \, |\tracefreeRic|^2 
&\geq \sigma \, |\Ric|^2 \, |\tracefreeRic|^2 \\ 
&\geq \frac{1}{3} \, \sigma \, \scal^2 \, |\tracefreeRic|^2 \\ 
&= \frac{1}{3} \, \sigma \, \scal^{4-\sigma} \, f \\ 
&\geq \frac{1}{3} \, \sigma \, \scal^{3-\sigma} \, f^{1+\frac{1}{\sigma}}. 
\end{align*} 
Putting these facts together, we conclude that 
\[\frac{\partial}{\partial t} f \leq \Delta f + \frac{2(1-\sigma)}{\scal} \, \partial_k \scal \, \partial^k f - \frac{2}{3} \, \sigma \, f^{1+\frac{1}{\sigma}}.\] 
Using the maximum principle, we obtain 
\[f \leq \Big ( \frac{3}{2t} \Big )^\sigma.\]
This completes the proof. \\

\begin{corollary}
Let $M$ be a compact three-manifold, and let $g(t)$, $t \in (-\infty,0)$, be an ancient solution to the Ricci flow on $M$. Moreover, suppose that there exists a uniform constant $\rho \in (0,1)$ such that 
\[\Ric_{g(t)} \geq \rho \, \scal_{g(t)} \, g(t) \geq 0\] 
for each $t \in (-\infty,0)$. Then the manifold $(M,g(t))$ has constant sectional curvature for each $t \in (-\infty,0)$. 
\end{corollary}

\textbf{Proof.} 
It follows from Proposition \ref{interior.estimate.dim.3} that $|\tracefreeRic_{g(t)}|^2 = 0$ for each $t \in (-\infty,0)$. Therefore, the manifold $(M,g(t))$ has constant sectional curvature for each $t \in (-\infty,0)$.

\section{The higher dimensional case}

\label{general.tools}

In this section, we develop some general tools that will be used in proof of Theorem \ref{ancient.solutions.higher.dim}. To that end, we fix an integer $n \geq 4$. Moreover, we denote by $\mathscr{C}_B(\mathbb{R}^n)$ the space of algebraic curvature tensors on $\mathbb{R}^n$. Given any algebraic curvature tensor $R \in \mathscr{C}_B(\mathbb{R}^n)$, we define an algebraic curvature tensor $Q(R) \in \mathscr{C}_B(\mathbb{R}^n)$ by 
\[Q(R)_{ijkl} = \sum_{p,q=1}^n R_{ijpq} \, R_{klpq} + 2 \sum_{p,q=1}^n (R_{ipkq} \, R_{jplq} - R_{iplq} \, R_{jpkq}).\] 
The expression $Q(R)$ arises naturally in the evolution equation for the curvature tensor under Ricci flow (cf. \cite{Hamilton2}; see also \cite{Brendle-book}, Section 2.3). The ordinary differential equation $\frac{d}{dt} R = Q(R)$ on the space $\mathscr{C}_B(\mathbb{R}^n)$ will be referred to as the Hamilton ODE.

We next consider a cone $C \subset \mathscr{C}_B(\mathbb{R}^n)$. We say that the cone $C$ has property $(*)$ if the following conditions are met:
\begin{itemize}
\item[(i)] $C$ is closed, convex, and $O(n)$-invariant.
\item[(ii)] $C$ is transversally invariant under the Hamilton ODE $\frac{d}{dt} R = Q(R)$.
\item[(iii)] Every algebraic curvature tensor $R \in C \setminus \{0\}$ has positive scalar curvature.
\item[(iv)] The curvature tensor $I_{ijkl} = \delta_{ik} \, \delta_{jl} - \delta_{il} \, \delta_{jk}$ lies in the interior of $C$.
\end{itemize}

In the remainder of this section, we assume that $C \subset \mathscr{C}_B(\mathbb{R}^n)$ is a cone satisfying $(*)$. Then $Q(R)$ lies in the interior of the tangent cone $T_R C$ for all $R \in C \setminus \{0\}$. By continuity, we can find a real number $\alpha_0 > 0$ such that 
\[Q(R + \alpha \, \scal(R) \, I) - \alpha_0^2 \, \scal(R)^2 \, I \in T_R C\] 
for all $R \in C \setminus \{0\}$ and all $\alpha \in [0,\alpha_0]$. Moreover, there exists a real number $\Lambda > 0$ such that $|\Ric(R)| \leq \Lambda \, \scal(R)$ for all $R \in C$. Let 
\[\delta = \min \Big \{ \frac{1}{2n(n-1)},\frac{\alpha_0}{2},\frac{\alpha_0^2}{4(1+2\Lambda^2)} \Big \} > 0.\] 
For each $t \in [0,\delta]$, we define a subset $F(t) \subset \mathscr{C}_B(\mathbb{R}^n)$ by 
\[F(t) = \{R \in C: R + (1 - t \, \scal(R)) \, I \in C\}.\] 
Clearly, $F(t)$ is closed, convex, and $O(n)$-invariant. Moreover, $F(0) = C$.

\begin{lemma}
\label{auxiliary.result}
Suppose that $R$ is an algebraic curvature tensor on $\mathbb{R}^n$ such that $R \in C$ and $R + (1 - t \, \scal(R)) \, I \in C$ for some $t \in [0,\delta]$. Then 
\[Q(R) - \scal(R) \, I - 2t \, |\Ric(R)|^2 \, I\] 
lies in the interior of the tangent cone to $C$ at the point $R + (1 - t \, \scal(R)) \, I$.
\end{lemma}

\textbf{Proof.} 
If $t \, \scal(R) < 1$, then the sum $R + (1 - t \, \scal(R)) \, I$ lies in the interior of $C$. In this case, the assertion is trivial.

Hence, it suffices to consider the case $t \, \scal(R) \geq 1$. For abbreviation, let 
\[S = R + (1-t \, \scal(R)) \, I \in C.\] 
Since $t \in [0,\delta]$, we have 
\[\scal(S) > (1 - n(n-1)t) \, \scal(R) \geq \frac{1}{2} \, \scal(R).\] 
Hence, if we put 
\[\alpha = \frac{t \, \scal(R) - 1}{\scal(S)},\] 
then we have $0 \leq \alpha < 2t \leq \alpha_0$. Since $S \in C \setminus \{0\}$, it follows that 
\[Q(S + \alpha \, \scal(S) \, I) - \alpha_0^2 \, \scal(S)^2 \, I \in T_S C\] 
by definition of $\alpha_0$. We next observe that 
\[S + \alpha \, \scal(S) \, I = R\] 
and 
\[\alpha_0^2 \, \scal(S)^2 > \frac{\alpha_0^2}{4} \, \scal(R)^2 \geq (1 + 2\Lambda^2) \, t \, \scal(R)^2 \geq \scal(R) + 2t \, |\Ric(R)|^2.\] 
Putting these facts together, we conclude that 
\[Q(R) - \scal(R) \, I - 2t \, |\Ric(R)|^2 \, I\] 
lies in the interior of the tangent cone $T_S C$. This completes the proof. \\

\begin{proposition}
\label{ode}
Suppose that $R(t)$ is a solution of the Hamilton ODE $\frac{d}{dt} R(t) = Q(R(t))$ which is defined on some time interval $[t_0,t_1] \subset [0,\delta]$. If $R(t_0) \in F(t_0)$, then $R(t) \in F(t)$ for all $t \in [t_0,t_1]$.
\end{proposition}

\textbf{Proof.} 
By assumption, we have $R(t_0) \in C$. Since $C$ is invariant under the Hamilton ODE, we conclude that $R(t) \in C$ for all $t \in [t_0,t_1]$. Hence, it suffices to show that $R(t) + (1 - t \, \scal(R(t))) \, I \in C$ for all $t \in [t_0,t_1]$. 

For abbreviation, let 
\[S(t) = R(t) + (1 - t \, \scal(R(t))) \, I\] 
for all $t \in [t_0,t_1]$. Since $R(t)$ is a solution of the Hamilton ODE, we have 
\[\frac{d}{dt} S(t) = Q(R(t)) - \scal(R(t)) \, I - 2t \, |\Ric(R(t))|^2 \, I\] 
for all $t \in [t_0,t_1]$. We claim that $S(t) \in C$ for all $t \in [t_0,t_1]$. Suppose this false. We define a real number $\tau$ by 
\[\tau = \inf \{t \in [t_0,t_1]: S(t) \notin C\}.\] 
By definition of $\tau$, we have $\tau \in [0,\delta]$ and $S(\tau) \in C$. Furthermore, we have $R(\tau) \in C$. Hence, Lemma \ref{auxiliary.result} implies that the derivative $\frac{d}{dt} S(t) \big |_{t=\tau}$ lies in the interior of the tangent cone $T_{S(\tau)} C$. By Proposition 5.4 in \cite{Brendle-book}, there exists a real number $\varepsilon > 0$ such that $S(t) \in C$ for all $t \in [\tau,\tau+\varepsilon)$. This contradicts the definition of $\tau$. \\

\begin{corollary} 
\label{interior.estimate.higher.dim}
Let $\delta$ be defined as above. Moreover, let $g(t)$, $t \in [0,\delta]$, be a solution to the Ricci flow on a compact $n$-dimensional manifold $M$. Finally, we assume that the curvature tensor of $(M,g(0))$ lies in the cone $C$ for all points $p \in M$. Then 
\[R_{g(t)} + (1 - t \, \scal_{g(t)}) \, I \in C\] 
for all points $(p,t) \in M \times [0,\delta]$.
\end{corollary}

\textbf{Proof.} 
By assumption, the curvature tensor of $(M,g(0))$ lies in the set $F(0)$ for all points $p \in M$. Using Proposition \ref{ode} and the maximum principle (cf. \cite{Chow-Lu}, Theorem 3), we conclude that the curvature tensor of $(M,g(t))$ lies in the set $F(t)$ for all points $(p,t) \in M \times [0,\delta]$. This proves the assertion. \\

\begin{corollary} 
\label{ancient}
Let $g(t)$, $t \in (-\infty,0)$, be an ancient solution to the Ricci flow on a compact $n$-dimensional manifold $M$. Moreover, suppose that the curvature tensor of $(M,g(t))$ lies in the cone $C$ for all $t \in (-\infty,0)$. Then 
\[R_{g(t)} - \delta \, \scal_{g(t)} \, I \in C\] 
for all points $(p,t) \in M \times (-\infty,0)$.
\end{corollary}

\textbf{Proof.} 
Fix a time $\tau \in (-\infty,0)$ and a real number $\sigma > 0$. We define a one-parameter family of metrics $\tilde{g}(t)$, $t \in [0,\delta]$, by 
\[\tilde{g}(t) = \sigma \, g \Big ( \frac{t-\delta}{\sigma} + \tau \Big ).\] 
Clearly, the metrics $\tilde{g}(t)$, $t \in [0,\delta]$, form a solution to the Ricci flow. By assumption, the curvature tensor of $(M,\tilde{g}(0))$ lies in the cone $C$ for all points $p \in M$. Hence, it follows from Corollary \ref{interior.estimate.higher.dim} that 
\[R_{\tilde{g}(\delta)} + (1 - \delta \, \scal_{\tilde{g}(\delta)}) \, I \in C\] 
for all points $p \in M$. This implies 
\[R_{g(\tau)} + (\sigma - \delta \, \scal_{g(\tau)}) \, I \in C\] 
for all points $p \in M$. Taking the limit as $\sigma \to 0$, we conclude that 
\[R_{g(\tau)} - \delta \, \scal_{g(\tau)} \, I \in C\] 
for all points $p \in M$. Since $\tau \in (-\infty,0)$ is arbitrary, the assertion follows. \\

\begin{theorem}
\label{general.principle}
Let $C(s)$, $s \in [0,1]$, be a family of cones in $\mathscr{C}_B(\mathbb{R}^n)$ satisfying property $(*)$. Moreover, suppose that the cones $C(s)$ vary continuously in $s$. Finally, let $g(t)$, $t \in (-\infty,0)$, be an ancient solution to the Ricci flow on a compact $n$-dimensional manifold $M$ such that $R_{g(t)} \in C(0)$ for all points $(p,t) \in M \times (-\infty,0)$. Then $R_{g(t)} \in C(1)$ for all $(p,t) \in M \times (-\infty,0)$.
\end{theorem}

\textbf{Proof.} 
Let $\mathscr{S}$ denote the set of all real numbers $s \in [0,1]$ with the property that $R_{g(t)} \in C(s)$ for all points $(p,t) \in M \times (-\infty,0)$. We claim that $\mathscr{S} = [0,1]$. 

Clearly, $\mathscr{S}$ is closed and non-empty. We next show that $\mathscr{S}$ is an open subset of $[0,1]$. To that end, we fix a real number $s_0 \in \mathscr{S}$. Then $R_{g(t)} \in C(s_0)$ for all points $(p,t) \in M \times (-\infty,0)$. By Corollary \ref{ancient}, there exists a real number $\delta > 0$ such that 
\[R_{g(t)} - \delta \, \scal_{g(t)} \, I \in C(s_0)\] 
for all points $(p,t) \in M \times (-\infty,0)$. Since the cones $C(s)$ vary continuously in $s$, there exists a real number $\varepsilon > 0$ such that $R_{g(t)} \in C(s)$ for all points $(p,t) \in M \times (-\infty,0)$ and all $s \in [s_0-\varepsilon,s_0+\varepsilon] \cap [0,1]$. Consequently, we have $[s_0-\varepsilon,s_0+\varepsilon] \cap [0,1] \subset \mathscr{S}$. This shows that $\mathscr{S}$ is an open subset of $[0,1]$. Thus, we conclude that $\mathscr{S} = [0,1]$, as claimed.

\section{Proof of Theorem \ref{ancient.solutions.higher.dim}}

We now describe the proof of Theorem \ref{ancient.solutions.higher.dim}. As in the previous section, we fix an integer $n \geq 4$. We denote by $\tilde{C}$ and $\hat{C}$ the cones introduced in \cite{Brendle1} and \cite{Brendle-Schoen1}. The cone $\tilde{C}$ consists of all algebraic curvature tensors $R \in \mathscr{C}_B(\mathbb{R}^n)$ satisfying
\begin{align*} 
&R(e_1,e_3,e_1,e_3) + \lambda^2 \, R(e_1,e_4,e_1,e_4) \\ 
&+ R(e_2,e_3,e_2,e_3) + \lambda^2 \, R(e_2,e_4,e_2,e_4) \\ 
&- 2\lambda \, R(e_1,e_2,e_3,e_4) \geq 0 
\end{align*} 
for all orthonormal four-frames $\{e_1,e_2,e_3,e_4\} \subset \mathbb{R}^n$ and all $\lambda \in [0,1]$. Similarly, the cone $\hat{C}$ consists of all algebraic curvature tensors $R \in \mathscr{C}_B(\mathbb{R}^n)$ satisfying
\begin{align*} 
&R(e_1,e_3,e_1,e_3) + \lambda^2 \, R(e_1,e_4,e_1,e_4) \\ 
&+ \mu^2 \, R(e_2,e_3,e_2,e_3) + \lambda^2 \mu^2 \, R(e_2,e_4,e_2,e_4) \\ 
&- 2\lambda \mu \, R(e_1,e_2,e_3,e_4) \geq 0 
\end{align*} 
for all orthonormal four-frames $\{e_1,e_2,e_3,e_4\} \subset \mathbb{R}^n$ and all $\lambda,\mu \in [0,1]$. The cones $\tilde{C}$ and $\hat{C}$ are both invariant under the Hamilton ODE $\frac{d}{dt} R = Q(R)$. A detailed discussion of these cones can be found in \cite{Brendle-book}, Chapter 7.

We next describe a family of invariant curvature cones interpolating between the cone $\tilde{C}$ and the cone $\hat{C}$. For each $s \in (0,\infty)$, we denote by $\tilde{C}(s)$ the set of all algebraic curvature tensors $R \in \mathscr{C}_B(\mathbb{R}^n)$ such that 
\begin{align*} 
&R(e_1,e_3,e_1,e_3) + \lambda^2 \, R(e_1,e_4,e_1,e_4) \\ 
&+ \mu^2 \, R(e_2,e_3,e_2,e_3) + \lambda^2 \mu^2 \, R(e_2,e_4,e_2,e_4) \\ 
&- 2\lambda \mu \, R(e_1,e_2,e_3,e_4) + \frac{1}{s} \, (1 - \lambda^2) \, (1 - \mu^2) \, \scal(R) \geq 0 
\end{align*} 
for all orthonormal four-frames $\{e_1,e_2,e_3,e_4\} \subset \mathbb{R}^n$ and all $\lambda,\mu \in [0,1]$. Clearly, $\tilde{C}(s)$ is a closed, convex cone, which is invariant under the natural action of $O(n)$. Moreover, we have $\hat{C} \subset \tilde{C}(s) \subset \tilde{C}$ for each $s \in (0,\infty)$. The following result is an immediate consequence of Proposition 10 in \cite{Brendle1}:

\begin{proposition}
\label{interpolating.cones}
For each $s \in (0,\infty)$, the cone $\tilde{C}(s)$ is invariant under the Hamilton ODE $\frac{d}{dt} R = Q(R)$. 
\end{proposition}

\textbf{Proof.} 
Let us fix a real number $s \in (0,\infty)$. Moreover, let $R(t)$, $t \in [0,T)$, be a solution of the Hamilton ODE such that $R(0) \in \tilde{C}(s)$. We claim that $R(t) \in \tilde{C}(s)$ for all $t \in [0,T)$. Without loss of generality, we may assume that $\scal(R(0)) = s$. This implies 
\begin{align*} 
&R(0)(e_1,e_3,e_1,e_3) + \lambda^2 \, R(0)(e_1,e_4,e_1,e_4) \\ 
&+ \mu^2 \, R(0)(e_2,e_3,e_2,e_3) + \lambda^2 \mu^2 \, R(0)(e_2,e_4,e_2,e_4) \\ 
&- 2\lambda \mu \, R(0)(e_1,e_2,e_3,e_4) + (1 - \lambda^2) \, (1 - \mu^2) \geq 0 
\end{align*} 
for all orthonormal four-frames $\{e_1,e_2,e_3,e_4\} \subset \mathbb{R}^n$ and all $\lambda,\mu \in [0,1]$. Hence, Proposition 10 in \cite{Brendle1} implies that 
\begin{align*} 
&R(t)(e_1,e_3,e_1,e_3) + \lambda^2 \, R(t)(e_1,e_4,e_1,e_4) \\ 
&+ \mu^2 \, R(t)(e_2,e_3,e_2,e_3) + \lambda^2 \mu^2 \, R(t)(e_2,e_4,e_2,e_4) \\ 
&- 2\lambda \mu \, R(t)(e_1,e_2,e_3,e_4) + (1 - \lambda^2) \, (1 - \mu^2) \geq 0 
\end{align*} 
for all orthonormal four-frames $\{e_1,e_2,e_3,e_4\} \subset \mathbb{R}^n$, all $\lambda,\mu \in [0,1]$, and all $t \in [0,T)$. Since $\scal(R(t)) \geq \scal(R(0)) = s$, we conclude that $R(t) \in \tilde{C}(s)$ for all $t \in [0,T)$. \\

After these preparations, we now present the proof of Theorem \ref{ancient.solutions.higher.dim}. 

\begin{theorem}
\label{step.1}
Assume that $g(t)$, $t \in (-\infty,0)$, is an ancient solution to the Ricci flow on a compact $n$-dimensional manifold $M$. Moreover, we assume that there exists a uniform constant $\rho > 0$ such that 
\[R_{g(t)} - \rho \, \scal_{g(t)} \, I \in \hat{C}\] 
for all points $(p,t) \in M \times (-\infty,0)$. Then the manifold $(M,g(t))$ has constant sectional curvature for each $t \in (-\infty,0)$.
\end{theorem} 

\textbf{Proof.} 
Consider the one-parameter family of cones $\hat{C}(s)$, $s \in (0,\infty)$, defined in \cite{Brendle-Schoen1}. It is shown in \cite{Brendle-Schoen1} that the cone $\hat{C}(s)$ has property $(*)$ for each $s \in (0,\infty)$. Furthermore, the cones $\hat{C}(s)$ vary continuously in $s$. 

By assumption, there exists a real number $s_0 \in (0,\infty)$ such that $R_{g(t)} \in \hat{C}(s_0)$ for all points $(p,t) \in M \times (-\infty,0)$. Using Theorem \ref{general.principle}, we conclude that $R_{g(t)} \in \hat{C}(s)$ for all points $(p,t) \in M \times (-\infty,0)$ and all $s \in (0,\infty)$. Consequently, the manifold $(M,g(t))$ has constant sectional curvature for each $t \in (-\infty,0)$. This completes the proof of Theorem \ref{step.1}. \\

\begin{theorem}
\label{step.2}
Assume that $g(t)$, $t \in (-\infty,0)$, is an ancient solution to the Ricci flow on a compact $n$-dimensional manifold $M$. Moreover, we assume that there exists a uniform constant $\rho > 0$ such that 
\[R_{g(t)} - \rho \, \scal_{g(t)} \, I \in \tilde{C}\] 
for all points $(p,t) \in M \times (-\infty,0)$. Then the manifold $(M,g(t))$ has constant sectional curvature for each $t \in (-\infty,0)$.
\end{theorem} 

\textbf{Proof.} 
By assumption, we have 
\[R_{g(t)} - \rho \, \scal_{g(t)} \, I \in \tilde{C}\] 
for all points $(p,t) \in M \times (-\infty,0)$. Hence, we can find a real number $s_0 \in (0,\infty)$ such that 
\[R_{g(t)} - \frac{1}{2} \, \rho \, \scal_{g(t)} \, I \in \tilde{C}(s_0)\] 
for all points $(p,t) \in M \times (-\infty,0)$.

We next consider a pair of real numbers $a,b$ such that $2a = 2b + (n-2)b^2$ and $b \in \big ( 0,\frac{\sqrt{2n(n-2)+4}-2}{n(n-2)} \big ]$. Following \cite{Bohm-Wilking}, we define a linear transformation $\ell_{a,b}: \mathscr{C}_B(\mathbb{R}^n) \to \mathscr{C}_B(\mathbb{R}^n)$ by 
\[\ell_{a,b}(R) = R + b \, \Ric(R) \owedge \text{\rm id} + \frac{1}{n} \, (a-b) \, \scal(R) \, \text{\rm id} \owedge \text{\rm id},\] 
where $\owedge$ denotes the Kulkarni-Nomizu product; see e.g. \cite{Besse}, Definition 1.110. If we choose $b \in \big ( 0,\frac{\sqrt{2n(n-2)+4}-2}{n(n-2)} \big ]$ sufficiently small, then 
\[R_{g(t)} \in \ell_{a,b}(\tilde{C}(s_0))\] 
for all points $(p,t) \in M \times (-\infty,0)$.

By Proposition \ref{interpolating.cones}, the cone $\tilde{C}(s)$ is invariant under the Hamilton ODE for each $s \in (0,\infty)$. Consequently, the cone $\ell_{a,b}(\tilde{C}(s))$ is transversally invariant under the Hamilton ODE for each $s \in (0,\infty)$ (cf. \cite{Bohm-Wilking}, Proposition 3.2). Therefore, the cone $\ell_{a,b}(\tilde{C}(s))$ has property $(*)$ for each $s \in (0,\infty)$. Using Theorem \ref{general.principle}, we conclude that $R_{g(t)} \in \ell_{a,b}(\tilde{C}(s))$ for all points $(p,t) \in M \times (-\infty,0)$ and all $s \in (0,\infty)$. Taking the limit as $s \to \infty$, we obtain $R_{g(t)} \in \ell_{a,b}(\hat{C})$ for all points $(p,t) \in M \times (-\infty,0)$. Hence, it follows from Theorem \ref{step.1} that $(M,g(t))$ has constant sectional curvature for each $t \in (-\infty,0)$.

\section{Ancient solutions satisfying a diameter bound}

In this final section, we study ancient solutions to the Ricci flow satisfying a suitable diameter bound. Throughout this section, we assume that $M$ is a compact manifold of dimension $n$, and $g(t)$, $t \in (-\infty,0)$, is a solution to the Ricci flow on $M$. The following proposition is a consequence of the differential Harnack inequality established in \cite{Brendle2}. 

\begin{lemma}
\label{corollary.of.harnack}
Suppose that the curvature tensor of $(M,g(t))$ lies in the cone $\hat{C}$ for each $t \in (-\infty,0)$. Then 
\begin{equation} 
\inf_M \scal_{g(\tau/2)} \geq \exp \Big ( -\frac{\text{\rm diam}(M,g(\tau))^2}{|\tau|} \Big ) \, \sup_M \scal_{g(\tau)}. 
\end{equation}
for all $\tau \in (-\infty,0)$.
\end{lemma}

\textbf{Proof.} 
Fix an arbitrary pair of points $p,q \in M$. We can find a smooth path $\gamma: [\tau,\tau/2] \to M$ such that $\gamma(\tau) = p$, $\gamma(\tau/2) = q$, and 
\[|\gamma'(t)|_{g(\tau)} = \frac{2 \, d_{g(\tau)}(p,q)}{|\tau|}.\] 
This implies 
\[|\gamma'(t)|_{g(t)} \leq \frac{2 \, d_{g(\tau)}(p,q)}{|\tau|}\] 
for all $t \in [\tau,\tau/2]$. Using the trace Harnack inequality (cf. \cite{Brendle2}, Proposition 13), we obtain 
\[\frac{\partial}{\partial t} \scal + 2 \, \partial_i \scal \, v^i \geq -2 \, \Ric(v,v)\] 
for every tangent vector $v$. Putting $v = \frac{1}{2} \, \gamma'(t)$ gives 
\begin{align*} 
\frac{d}{dt} \scal_{g(t)}(\gamma(t)) 
&\geq -\frac{1}{2} \, \Ric_{g(t)}(\gamma'(t),\gamma'(t)) \\ 
&\geq -\frac{1}{2} \, \scal_{g(t)}(\gamma(t)) \, |\gamma'(t)|_{g(t)}^2 \\ 
&\geq -\frac{2 \, d_{g(\tau)}(p,q)^2}{|\tau|^2} \, \scal_{g(t)}(\gamma(t)) 
\end{align*}
for all $t \in [\tau,\tau/2]$. Thus, we conclude that 
\[\scal_{g(\tau/2)}(q) \geq \exp \Big ( -\frac{d_{g(\tau)}(p,q)^2}{|\tau|} \Big ) \, \scal_{g(\tau)}(p).\] 
Since $p,q \in M$ are arbitrary, the assertion follows. \\

\begin{proposition}
\label{curvature.bound}
Suppose that the curvature tensor of $(M,g(t))$ lies in the cone $\hat{C}$ for each $t \in (-\infty,0)$. Moreover, suppose that 
\[\limsup_{\tau \to -\infty} \frac{1}{\sqrt{|\tau|}} \, \text{\rm diam}(M,g(\tau)) < \infty.\] 
Then 
\[\limsup_{\tau \to -\infty} \Big [ |\tau| \, \sup_M \scal_{g(\tau)} \Big ] < \infty.\] 
\end{proposition}

\textbf{Proof.} 
Since the solution $g(t)$ is defined until time $0$, we have 
\[\inf_M \scal_{g(\tau/2)} \leq \frac{n}{|\tau|}\] 
for each $\tau \in (-\infty,0)$ (see e.g. \cite{Brendle-book}, Proposition 2.19). Using Lemma \ref{corollary.of.harnack}, we deduce that 
\[\sup_M \scal_{g(\tau)} \leq \frac{n}{|\tau|} \, \exp \Big ( \frac{\text{\rm diam}(M,g(\tau))^2}{|\tau|} \Big ).\] From this, the assertion follows. \\

Finally, we recall the following result due to B.~Kostant (see \cite{Kostant}, Corollary 2.2):

\begin{proposition}
\label{holonomy}
Let $(N,h)$ be a compact, simply connected Riemannian manifold of dimension $n \neq 5$ which is, topologically, a rational homology sphere. Then the holonomy representation of $(N,h)$ is complete; that is, $(N,h)$ has holonomy group $\text{\rm SO}(n)$.
\end{proposition}

\begin{theorem}
Let $g(t)$, $t \in (-\infty,0)$, be an ancient solution to the Ricci flow on a compact, even-dimensional manifold $M$. Suppose that the curvature tensor of $(M,g(t))$ lies in the interior of the cone $\hat{C}$ for each $t \in (-\infty,0)$. Moreover, suppose that 
\[\limsup_{\tau \to -\infty} \frac{1}{\sqrt{|\tau|}} \, \text{\rm diam}(M,g(\tau)) < \infty.\] 
Then $(M,g(t))$ has constant sectional curvature for each $t \in (-\infty,0)$.
\end{theorem}

\textbf{Proof.} 
Suppose the assertion is false. By Theorem \ref{step.1}, we can find a sequence of points $(p_k,\tau_k) \in M \times (-\infty,0)$ such that $\lim_{k \to \infty} \tau_k = -\infty$ and 
\begin{equation} 
\label{curvature.approaching.boundary.of.C.hat}
R_{g(\tau_k)} - \frac{1}{k} \, \scal_{g(\tau_k)} \, I \notin \hat{C} 
\end{equation}
at $p_k$. For each $k$, we consider the rescaled metrics 
\[\tilde{g}_k(t) = \frac{1}{|\tau_k|} \, g(|\tau_k| \, t), \qquad t \in (-2,-\frac{1}{2}).\] 
For each $k$, the metrics $\tilde{g}_k(t)$, $t \in (-2,-\frac{1}{2})$, form a solution to the Ricci flow on $M$. By assumption, the diameter of $(M,\tilde{g}_k(t))$ has uniformly bounded diameter; moreover, it has uniformly bounded curvature by Proposition \ref{curvature.bound}. Since $M$ is even-dimensional, we conclude that the injectivity radius of $(M,\tilde{g}_k(t))$ is uniformly bounded from below.

Hence, after passing to a subsequence if necessary, the sequence $(M,\tilde{g}_k(t))$ converges in the Cheeger-Gromov sense to some limiting solution $(M,\bar{g}(t))$ to the Ricci flow. This limiting solution is defined for all $t \in (-2,-\frac{1}{2})$. Clearly, the curvature tensor of $(M,\bar{g}(t))$ lies in the cone $\hat{C}$ for each $t \in (-2,-\frac{1}{2})$. Moreover, it follows from (\ref{curvature.approaching.boundary.of.C.hat}) that the curvature tensor of $(M,\bar{g}(-1))$ lies on the boundary of the cone $\hat{C}$ for some point $q \in M$. By Proposition 9 in \cite{Brendle-Schoen2}, the manifold $(M,\bar{g}(-1))$ has non-generic holonomy group, i.e. $\text{\rm Hol}^0(M,\bar{g}(-1)) \neq \text{\rm SO}(n)$. On the other hand, it follows from the Differentiable Sphere Theorem that the universal cover of $(M,\bar{g}(-1))$ is diffeomorphic to $S^n$ (cf. \cite{Brendle-Schoen1}, Theorem 3). By Proposition \ref{holonomy}, the universal cover of $(M,\bar{g}(-1))$ has holonomy group $\text{\rm SO}(n)$. This is a contradiction.

\end{document}